\documentclass{amsart}
\usepackage{amsopn}
\usepackage{amsfonts}
\usepackage{hhline}

\newtheorem{thm}{Theorem}[subsection]
\newtheorem{cor}[thm]{Corollary}
\newtheorem{lem}[thm]{Lemma}
\theoremstyle{definition}
\newtheorem{defn}[thm]{Definition}
\theoremstyle{remark}
\newtheorem{rem}[thm]{Remark}
\DeclareMathOperator{\Lat}{Lat}

\begin{document}

\title
 {Principal Ideals in Subalgebras of Groupoid C*-algebras}

\author{Srilal Krishnan }

\address{Department of Mathematics,University of Alabama, Tuscaloosa}

\curraddr{P.O. Box 862547, Tuscaloosa, Al 35486-0024}

\email{krish001@bama.ua.edu}


\thanks{Thanks are due to my advisor Dr. Alan Hopenwasser}


\date{April 15, 2001.}

%%% ----------------------------------------------------------------------

\begin{abstract}
  The study of different types of ideals in non self-adjoint
   operator algebras has been a topic of recent research. This
   study focuses on principal ideals in subalgebras of groupoid
   C*-algebras. An ideal is said to be principal if it is
   generated by a single element of the algebra. We look at
   subalgebras of r-discrete principal groupoid C*-algebras and
   prove that these algebras are principal ideal algebras. Regular
   canonical subalgebras of almost finite C*-algebras have digraph
   algebras as their building blocks. The spectrum of almost
   finite C*-algebras has the structure of an r-discrete principal
   groupoid and this helps in the coordinization of these
   algebras. Regular canonical subalgebras of almost finite
   C*-algebras have representations in terms of open subsets of
   the spectrum for the enveloping C*-algebra. We conclude that
   regular canonical subalgebras are principal ideal algebras.

\end{abstract}

%%% ----------------------------------------------------------------------
\maketitle
% ------------------------------------------------------------------------
\section*{Introduction}
 Non self-adjoint limit algebras are direct limits of subalgebras
     of finite dimensional C*-algebras. The refinement embedding algebra,
     standard embedding algebra and the alternating embedding algebra are examples
     of non self-adjoint limit algebras and these examples are explained in detail in
     this paper. Some important special classes of non self-adjoint limit algebras are
     :\begin{itemize}
     \item TUHF(Triangular Uniformly Hyperfinite) algebras.
     \item TAF(Triangular approximately finite dimensional)
     algebras.
     \item Regular canonical subalgebras. \end{itemize}
     TUHF algebras form a subclass of TAF algebras and TAF algebras
     form a subclass of regular canonical subalgebras.
     The upper triangular
     complex matrix algebra $T_n$ is a basic building block for TUHF algebras whereas
     direct sums of upper triangular matrix algebras are the basic building blocks for
     the TAF algebras. Digraph algebras are the building blocks for
     regular canonical subalgebras of AF C*-algebras. The refinement embedding algebra,
     standard embedding algebra and the alternating embedding algebra are examples
     of TUHF limit algebras. The study of different types of ideals in non self-adjoint
    operator algebras has been a topic of recent research.
     We will study one of the basic types of ideals:\ a principal ideal
   in some non self-adjoint limit algebras. An ideal
   is said to be principal if it is generated
   by a single element of the algebra.
   In this study we will first analyze the structure of ideals in digraph algebras
   and prove that digraph algebras are principal ideal algebras. Since regular canonical subalgebras of AF C*-algebras
    are infinite dimensional analogues of digraph algebras, it is
    natural to expect these algebras to be principal ideal algebras. But it is observed that the
    proof does not follow naturally and we have to resort to the spectrum of
    these algebras in the proof. The spectrum of AF C*-algebras
    has the structure of an r-discrete principal groupoid and it
    is this groupoid and substructures of the groupoid which help
    in the coordinization of these algebras. Any closed subalgebra
    of an AF C*-algebra which contains the diagonal has a
    functional representation in terms of an open subset of the
    spectrum of the enveloping C*-algebra and we use this
    representation to prove that regular canonical subalgebras of AF C*-algebras are
    principal ideal algebras.
%-----------------------------------------------------------------------------------

\section{Preliminaries}
     In this section we define AF C*-algebras, UHF C*-algebras, TAF algebras, TUHF
     algebras and regular canonical subalgebras of AF C*-algebras.
     TAF and TUHF algebras are non self-adjoint versions of AF C*-algebras and UHF C*-algebras
     respectively.
     We will illustrate the definitions with some examples. Also
     there are two equivalent ways of defining these algebras,
     one as an inductive limit and other as the closure of an
     increasing union of finite dimensional algebras and these are
     explained below.
     \bigskip
     \goodbreak
     \subsection{Direct limit of C*-algebras}
     \begin{defn} A C*-algebra is a
     norm closed self-adjoint subalgebra of the operator algebra B(H), for some Hilbert space H.\end{defn}
      Every finite dimensional C*-algebra is *-isomorphic to the direct sum of
     full matrix algebras (A proof of this is indicated in ~\cite{krd96}, page 74).
     \smallskip
      Next let ${B_1}\stackrel{\varphi_{1}}{\rightarrow} {B_2}\stackrel{\varphi_{2}}
     {\rightarrow}{B_3}\stackrel{\varphi_{3}}{\rightarrow}{
     B_4}\stackrel{\varphi_{4}}{\rightarrow} \cdots {B_n}\stackrel{\varphi_{n}}{\rightarrow}{B_{n+1}}
     \cdots$ denote an injective direct system of C*-algebras
     $B_n,n=1,2,\cdots$ with star injections $\varphi_{n}:B_{n}
     \rightarrow B_{n+1}$. Then the product $\prod_{n=1}^{\infty} B_n$
     is a *-algebra with pointwise defined operations. Let
     $B_{\infty}^{0}= \{b:b = (b_n) \in \prod_{n=1}^{\infty} B_n,
     \varphi_{n}(b_n)=b_{n+1}$, for all large $n$\}. Then
     $B_{\infty}^{0}$ is a *-subalgebra of $\prod_{n=1}^{\infty}
     B_n$. Any *-homomorphism $\varphi_{n}:B_{n}
     \rightarrow B_{n+1}$ is necessarily norm-decreasing (this is a
     standard result from C*-algebras). Consequently $\|b_{n+1}\| \leq \|b_{n}\|$ for
     large $n$ and so the sequence $(\|b_{n}\|)$ is eventually
     decreasing and also bounded below. Thus the sequence
     $(\|b_{n}\|)$ converges. Let $p(b)= \lim_{n \rightarrow
     \infty} \|b_{n}\|$. Then it is clear that $p:B_{\infty}^{0}
     \rightarrow R^{+}$, $b \mapsto p(b)$, is a C*-seminorm on
     $B_{\infty}^{0}$. We denote the enveloping C*-algebra of
     $(B_{\infty}^{0},p)$ by $B$, and call it the $direct$ $limit$
     of the sequence $(B_n,\varphi_{n})_{n=1}^{\infty}$. Also we
     denote $B$ by ${\displaystyle\lim_{n}{B_n}}$.
     \bigskip
     \goodbreak
     Next we define AF C*-algebra and UHF C*-algebra.
     \begin{defn} A C*-algebra $B$
     is approximately finite
     dimensional (AF) if it is the closure of an increasing union of finite-dimensional
     C*-subalgebras $B_n$.\end{defn} \begin{rem}Equivalently every AF C*-algebra is *-isomorphic
     to a direct limit algebra,  ${\displaystyle\lim_{i}(B_{n_i},\phi_{i})}$, associated with
     standard unital injective maps $\phi_{i}:{B_{n_i}}{\rightarrow}\
     {B_{n_{i+1}}}$, where each ${B_{n_i}}$ is a direct sum of full
     matrix algebras.\end{rem} A simple finite dimensional C*-algebra is *-isomorphic to
     a full matrix algebra. \begin{defn}A C*-algebra $B$ is uniformly hyperfinite
    (UHF) if it is the closure of an increasing union of simple
     finite-dimensional C*-subalgebras $B_n$.\end{defn} \begin{rem}Equivalently, every UHF algebra  is *-isomorphic
     to a limit algebra,  ${\displaystyle\lim_{i}(M_{n_i},\phi_{i})}$, associated with
     standard unital injective maps $\phi_{i}:{M_{n_i}}{\rightarrow}\ {M_{n_{i+1}}}$.\end{rem}
     \bigskip
     \goodbreak
     \subsection{Matrix unit system for limit algebras}
     Since  ${M_n}\cong{B(\mathbb{C}^n)}$, \ where ${\mathbb{C}^n}$ denotes the $n$-dimensional complex
     field,\ we will make frequent use of the fact that corresponding to any orthonormal
     basis ${e_{1},e_{2},\cdots,e_{n}}$,\ ${M_n}$ has a  basis consisting of
     matrix units ${E_{ij}}= e_{i}e_{j}^*$ ($e_j^*$ denotes the conjugate transpose of $e_j$)
     for $1 \le i,j \le n$. We can also define $e_{i}e_{j}^*$ as an operator:
     $e_{i}e_{j}^*(x)=<x,e_{j}>e_{i}$.
      Consequently if $B \cong {M_{n_1}}\oplus{M_{n_2}}\oplus{M_{n_3}}\oplus\cdots{M_{n_k}}$,
     then $B$ has a matrix unit system(m.u.s), say
     $\{E_{ij}^{s}:1 \le s \le k, 1 \le i,j \le n_{s}\}$. Evidently such a m.u.s is not unique because if $\{E_{ij}\}$
     is a m.u.s then $\{e^{{i\theta}(j-i)}E_{ij}\}$ is another m.u.s where $e^{i\theta}$ denotes a complex number of
     modulus 1. Also conjugating the m.u.s by any unitary will yield other matrix unit systems. But any two m.u.s of a finite
     dimensional C*-algebra are inner conjugate and hence the choice of m.u.s for a chain of finite
     dimensional C*-algebras is irrelevant although it does matter how the m.u.s fits with the embeddings.
     The m.u.s for a chain $\{B_k\}$ is a system
     $\{e_{ij}^k\}_{ijk}$ where for each $k$ the system $\{e_{ij}^k\}_{ij}$ is a m.u.s for $B_k$ and
     where each  ${e_{ij}^k}$ is a sum of the elements of $\{e_{ij}^{k+1}\}_{ij}$.
     Suppose $B=\overline{\cup_{n=1}^\infty {A_n}}$ is an AF algebra
     and let $D_n$ be a maximal abelian self-adjoint subalgebra (masa) for each $n$ such that
     $D_n \subseteq D_{n+1}$, for all $n$.
     Let $D=\overline{\cup_{i=1}^\infty {D_n}}$.
     Then $D$ is a masa of $B$ and $D_n=D\cap{B_n}$ is such that $D_n \subseteq D_{n+1}$. The existence
     of such a masa is guaranteed in ~\cite{jg60}. Let
     $B_n=\oplus_{m=1}^{l(n)}M_{k(n,m)}$, where $M_{k}$ denotes a $k \times k$ matrix.
     Then for each $n$ and $m$,\ a
     m.u.s $\{e_{ij}^{nm}\}$ can always be chosen for $M_{k(n,m)}$
     so that if $\phi_{n}:{B_n}{\rightarrow}\ {B_{n+1}}$
     denotes the embedding from ${B_n}$ to ${B_{n+1}}$
     then  $\phi_{n}({e_{ij}^{(nm)}})$ is a sum of matrix units of
     ${B_{n+1}}$. Also the m.u.s can be chosen such that each $D_n$ is generated by the diagonal matrix units.
     Consequently $D$ is the closed linear span of
     $\{e_{ii}^{(nm)}:1 \le n, 1\le m \le l(n), 1 \le i \le
     k(n,m)\}$.  (The reader is referred to ~\cite{scp92bk} for
     details.) \begin{rem} All subalgebras of AF algebras in this
     paper are norm closed.\end{rem} Let $D$ be the abelian C*-algebra generated
     by all of the diagonal matrix units $e_{ii}^{(nm)}$ associated with the m.u.s
     $\{e_{ii}^{(nm)}\}$, as above. Then $D=\overline{\cup_{i=1}^\infty
     {D_n}}$, where $D_{n}$ is the masa in $B_n$ spanned by the
     diagonal matrix units. Then $D$ is a masa in $B$ (for a proof of this, refer to
     ~\cite{scp92bk}). We call $D$ a $regular$ $canonical$ $masa$ associated
     with the m.u.s. $\{e_{ij}^{(nm)}\}$.
     \bigskip
     \goodbreak
     At this point we will define
     TAF algebras and TUHF algebras and regular canonical subalgebras.
     Let $A$ denote a regular canonical subalgebra of an AF C*-algebra $B$.  It is important to note that the
     embeddings $\phi_{i}:{A_{n_i}}{\rightarrow}\
     {A_{n_{i+1}}}$ are *-extendible and maps the normaliser of
     $D_k$ into the normaliser $D_{k+1}$. Such embeddings are
     called $regular$ $embeddings$ in literature.
     \begin{defn} A regular canonical subalgebra $A$ of an AF C*-algebra $B$ is a closed subalgebra of $B$
     such that $ D \subseteq A \subseteq B$;  where $D$ is a regular canonical masa associated with a m.u.s for $B$.
     \end{defn}
      \begin{defn} If $B=\overline{\cup_{i=1}^\infty
     {B_n}}$ is an AF algebra with masa $D$ , then a subalgebra $A$
     of $B$ is said to be TAF with diagonal $D$ if $D=A \cap A^*$.\end{defn}
      \begin{defn} If $B=\overline{\cup_{i=1}^\infty
     {B_n}}$ is a UHF algebra with masa $D$ , then a subalgebra $A$
     of $B$ is said to be TUHF with diagonal $D$ if $D=A \cap A^*$.\end{defn}
     \begin{defn} A triangular subalgebra $A$ of an AF C*-algebra
     $B$ is a closed subalgebra of $B$ such that $A \cap A^*$ is a
     masa.\end{defn}
     \begin{defn} A TAF subalgebra $A$ of an AF C*-algebra $B$ is said to be
     maximal triangular if $A$ is the only triangular subalgebra of $B$
     containing $A$.\end{defn} \begin{rem} In the above
     definition, if the sequence $\{B_n\}$ can be chosen such that $A \cap
     B_n$ is maximal triangular in $B_n$ for each $n$, then $A$ is
     called strongly maximal triangular.\end{rem} The refinement embedding algebra,
     standard embedding algebra and the alternating embedding algebra are
     examples of strongly maximal triangular algebras. It is evident that a strongly
     maximal TAF algebra is a maximal TAF algebra. But the converse is not true and
     an example is given in ~\cite{ppw90}.\begin{rem} If $A$
     is a TUHF algebra then it may be possible to write
     $A=\overline{\cup_{i=1}^\infty {A_n}}$ where each $A_n$ is not a factor.
     This motivates the following definition. \end{rem} \begin{defn} A
     strongly maximal triangular subalgebra $A$ of a UHF algebra
     $B$ is said to be strongly maximal in factors if a sequence $\{B_n\}$ can be
     chosen such that $B_n \cong {M_{n_k}}$ for each $n$, $B=\overline{\cup_{i=1}^\infty {B_n}}$
     and $A \cap B_n$ is maximal triangular in $B_n$ for each $n$.
     \end{defn} Again it is not true in general that a strongly
     maximal TUHF algebra is strongly maximal in factors. An
     example is given in ~\cite{ppw90}.
     \bigskip
     \goodbreak
     \subsection{Examples of limit algebras}
     We will study the refinement embedding algebra,
     standard embedding algebra, the alternating embedding algebra and digraph algebras in detail.
    \begin{enumerate}
     \item Let $(n_k)$ denote a sequence of positive integers such
     that $n_k$ divides $n_{k+1}$, for each $k=1,2,\cdots $ Consider the unital
     injective maps $\rho_{k}:{M_{n_k}}{\rightarrow}\
     {M_{n_{k+1}}}$ given by $\rho_{k}(a_{ij})=(a_{ij}I_{r_k})$
     such that $(a_{ij}I_{r_k})$ is the partitioned matrix in
     ${M_{n_{k+1}}}$, with $I_{r_k}$ the identity matrix in
     $M_{r_k}$, where $r_k=\frac{n_{k+1}}{n_k}$. In this case
     $\rho_{k}(e_{ij}^k)= \sum_{t=1}^{q_k}e_{(i-1){q_k}+t,(j-1){q_k}+t}^{k+1}$.
     Here is an example of a  refinement embedding :

     $\rho_{k}(\left[\begin{array}{cc}
        a & b \\
        c & d \\
   \end{array}\right])=\left[\begin{tabular}{cccccc|cccccc}
        a & 0 & . & . & 0 & 0 & b & 0 & . & . & 0 & 0 \\
        0 & a & . & . & 0 & 0 & 0 & b & . & . & 0 & 0\\
        . & . & . & . & 0 & 0 & . & . & . & . & 0 & 0\\
        . & . & . & . & 0 & 0 & . & . & . & . & 0 & 0\\
        0 & 0 & 0 & 0 & a & 0 & 0 & 0 & 0 & 0 & b & 0\\
        0 & 0 & 0 & 0 & 0 & a & 0 & 0 & 0 & 0 & 0 & b\\ \hhline{------|------}
        c & 0 & . & . & 0 & 0 & d & 0 & . & . & 0 & 0\\
        0 & c & . & . & 0 & 0 & 0 & d & . & . & 0 & 0\\
        . & . & . & . & 0 & 0 & . & . & . & . & 0 & 0\\
        . & . & . & . & 0 & 0 & . & . & . & . & 0 & 0\\
        0 & 0 & 0 & 0 & c & 0 & 0 & 0 & 0 & 0 & d & 0\\
        0 & 0 & 0 & 0 & 0 & c & 0 & 0 & 0 & 0 & 0 & d\\
     \end{tabular}\right]$

     where the righthand side is a $2$ by $2$ matrix, with each
     entry amplified by order ${r}$ using the identity matrix $I_{r}$.
     Then the limit algebra ${\displaystyle\lim_{k}(M_{n_k},\rho_{k})}$ associated with
     these unital injective maps $\rho_{k}:{M_{n_k}}{\rightarrow}\
     {M_{n_{k+1}}}$ is an example of a UHF C*-algebra (also called Glimm Algebra in literature).
     Then if $T_n \subseteq M_n$ denotes the algebra of upper triangular
     matrices relative to the standard m.u.s it follows that
     $\rho_{k}(T_{n_k})\subseteq T_{n_{k+1}}$. The canonical subalgebra
     ${\displaystyle\lim_{k}(T_{n_k},\rho_{k})}$ of the UHF
     C*-algebra is called a refinement limit algebra.
     \item Let $(n_k)$ denote a sequence of positive integers such
     that $n_k$ divides $n_{k+1}$, for each $k=1,2,\cdots$ Consider the unital
     injective maps $\sigma_{k}:{M_{n_k}}{\rightarrow}\
     {M_{n_{k+1}}}$ given by $\sigma_{k}(e_{ij}^k)=
     \sum_{t=0}^{q_k-1}e_{i+tq_k,j+tq_k}^{k+1}$. The unital injective maps $\sigma_{k}:{M_{n_k}}{\rightarrow}\
     {M_{n_{k+1}}}$ are called standard embeddings. Here is an
     example of a standard embedding :

    $\sigma_{k}(\left[\begin{array}{cc}
        a & b \\
        c & d \\
 \end{array}\right])=\left[\begin{array}{cccccccccccc}
        a & b & 0 & 0 & 0 & 0 & . & . & 0 & 0 & 0 & 0\\
        c & d & 0 & 0 & 0 & 0 & . & . & 0 & 0 & 0 & 0\\
        0 & 0 & a & b & 0 & 0 & . & . & 0 & 0 & 0 & 0\\
        0 & 0 & c & d & 0 & 0 & . & . & 0 & 0 & 0 & 0\\
        0 & 0 & 0 & 0 & a & b & . & . & 0 & 0 & 0 & 0\\
        0 & 0 & 0 & 0 & c & d & . & . & 0 & 0 & 0 & 0\\
        . & . & . & . & . & . & . & . & . & . & . & .\\
        . & . & . & . & . & . & . & . & . & . & . & .\\
        0 & 0 & 0 & 0 & 0 & 0 & . & . & a & b & 0 & 0\\
        0 & 0 & 0 & 0 & 0 & 0 & . & . & c & d & 0 & 0\\
        0 & 0 & 0 & 0 & 0 & 0 & . & . & 0 & 0 & a & b\\
        0 & 0 & 0 & 0 & 0 & 0 & . & . & 0 & 0 & c & d\\
     \end{array}\right]$

     Then if $T_n \subseteq M_n$ denotes the algebra of upper triangular
     matrices relative to the standard m.u.s it follows that
     $\sigma_{k}(T_{n_k})\subseteq T_{n_{k+1}}$. The canonical subalgebra
     ${\displaystyle\lim_{k}(T_{n_k},\rho_{k})}$ is called a standard limit algebra. \item
     Alternation limit algebras are limit algebras in which the
     refinement embeddings and the standard embeddings are used
     alternatively. Let $(r_k)$ and $(s_k)$ denote
     sequences of positive integers; $k=1,2,\cdots$. Let $\sigma_{k}:{M_{s_k}}{\rightarrow}\
     {M_{s_{k+1}}}$  denote the standard embeddings and $\rho_{k}:{M_{r_k}}{\rightarrow}\
     {M_{r_{k+1}}}$ denote the refinement embeddings. Then the
     limit of the direct system $C \stackrel{\rho}{\rightarrow}
     T_{r_1}\stackrel{\sigma}{\rightarrow}
     T_{r_1s_1}\stackrel{\rho}{\rightarrow}
     T_{r_1s_1r_2}\stackrel{\sigma}{\rightarrow}\cdots$ is called
     an alternation limit algebra. \item A digraph algebra is a
     subalgebra of the full complex matrix algebra $M_n$ which
     contains a maximal abelian subalgebra of $M_n$. Let $D_n$ denote
     the standard diagonal algebra associated with the standard m.u.s .
     The digraph algebra is unitarily equivalent to an algebra containing
     $D_n$. Thus an example of a digraph algebra is

     $\left[\begin{array}{ccccccc}
     * & 0 & 0 & 0 & 0 & 0 & 0 \\
     0 & * & * & * & * & * & 0 \\
     0 & * & * & * & 0 & * & 0 \\
     0 & * & * & * & 0 & * & 0 \\
     0 & * & 0 & 0 & * & * & 0 \\
     0 & * & * & * & * & * & 0 \\
     0 & 0 & 0 & 0 & 0 & 0 & * \
   \end{array}\right]$.

     Digraph algebras are basic building blocks for
     regular canonical subalgebras. Let $A$ denote  a digraph
     algebra such that $A \subseteq M_n$. Then we have
     $D_n \subseteq A \subseteq M_n$, where $D_n$ denotes
     the algebra of diagonal matrices. Let $\Lat A$ denote the
     lattice of invariant projections of $A$. Then $\Lat A= \{p \in
     M_n: p = p^{*}=p^2$ and $ap=pap,\forall a \in A\}$.
     Let $p \in \Lat A$. Then $p$ is also invariant for $D_n$. But $D_n$ is a
     self-adjoint algebra and so $p$ is invariant for $D_n$ implies
     $p$ lies in the $commutant$ of $D_n$, where $commutant$ is the
     set of all elements of $M_n$ that commute with all elements
     of $D_n$ and is denoted by $D_n^{c}$. But $D_n$ is a masa and
     so $D_n^{c}=D_n$. So $p \in D_n$. Since $p$ is an arbitrary element
     in $\Lat A$, we can infer that $\Lat A \subseteq D_n$ and consequently is a
     commutative lattice. Note that $\Lat A$ is a commutative lattice with
     $p \vee q = p + q -pq$ and
     $p \wedge q = pq$ for all elements $p,q \in \Lat A$. Thus the lattice of invariant
     projections for a digraph algebra is a commutative subspace
     lattice(CSL). So the digraph algebras are just the finite
     dimensional CSL algebras.
     \end{enumerate}
  \begin{rem} If $\frac {n_{k+1}}{n_{k}}=2$ for $k = 1,2,\cdots$
  in examples 1 and 2 above, then the standard limit algebra and
  the refinement limit algebra are called $2^\infty$ TUHF
  algebras.\end{rem}
      \bigskip
     \goodbreak
     \section{Isometric isomorphism of limit algebras}
     The UHF algebras corresponding to the refinement and standard
     embeddings are isometrically isomorphic. This is proved in
     ~\cite{jg60} by J.Glimm. In fact Glimm proved that the isomorphism class was independent
     of the nature of the embeddings and depended only on the dimensions of the finite dimensional
     factors. But this is not true in the case of  TUHF
     algebras. In fact the standard limit algebra and the
     refinement limit algebra are not isometrically isomorphic.
     Let $S$ and $R$ denote the $2^\infty$ TUHF
     algebras via the standard embedding and refinement embedding respectively
     and let $\Lat S$ and $\Lat R$ denote the lattice of invariant
     projections of $S$ and $R$ respectively. Then
     $\Lat S =\{\mathrm{0,1}\}$ and $ \Lat R =\mathrm{L}$ where
     $  L$= $\{\sum_{i=1}^{j}e_{ii}^n :1 \le j \le 2^{n}, 1 \le n
     \le \infty\}$ $ \cup \{0\}$, a nest. Thus the limit algebras are not
     isomorphic, since  $\Lat S$ and $\Lat R$ are not
     isomorphic. Also it is interesting to observe that the
     refinement embedding maps $\Lat T_{n_k}$ into $\Lat
     T_{n_{k+1}}$ whereas the standard embedding does not. The
     classification of refinement limit algebra,
     standard limit algebra and the alternating limit algebra up
     to isometric isomorphism has been done in ~\cite{hp92}. For all
     classifications the authors used the spectrum, also called
     the topologized fundamental relation. The results in this
     paper also make use of the spectrum. The spectrum of these
     algebras has been studied in detail in ~\cite{scp92bk}. But for the sake
     of completeness we will describe the spectrum of these
     algebras.
%%% ----------------------------------------------------------------------
 \bigskip
 \goodbreak
\section{Spectrum or Fundamental relation}

%% ----------------------------------------------------------------------
    \subsection{Introduction}
     To understand the concept of spectrum for a limit algebra, we
     will start by defining a normalising partial isometry and we
     will see how normalising partial isometries act on the maximal
     ideal space of the canonical masa of a limit algebra to
     generate the spectrum.
     \bigskip
     \goodbreak
     \subsection{Spectrum of limit algebras}
     Let $B$ denote a limit algebra ${\displaystyle\lim_{i}(B_i)}$
     arising from the direct system ${B_1}\stackrel{\varphi_{1}}{\rightarrow} {B_2}\stackrel{\varphi_{2}}
     {\rightarrow}{B_3}\stackrel{\varphi_{3}}{\rightarrow}{
     B_4}\stackrel{\varphi_{4}}{\rightarrow} \cdots {B_n}\stackrel{\varphi_{n}}{\rightarrow}{B_{n+1}}
     \cdots$ such that the *-extendible embeddings
     $\varphi_{i}$ map the matrix units of $B_i$ to sums of matrix units of
     $B_{i+1}$. Let $D$ denote the regular canonical masa of $B$
     associated with a m.u.s. Then with $D_n=D\cap{B_n}$, $D_n \subseteq
     D_{n+1}$ and $D=\overline{\cup_{n=1}^\infty {D_n}}$.
     \begin{defn} An element $p$ in $B$ is a projection if
     $p^*=p=p^2$. \end{defn} \begin{defn}An element $v$ in $B$ is a
     partial isometry if  $v^{*}v$ is a projection.\end{defn}
     \begin{defn} The range (final) projection and the
    domain (initial) projection of a partial isometry $v$ in $B$ is
    defined as $r(v)=vv^*$ and $d(v)=v^*v$ respectively.\end{defn} \begin{defn}A map $\alpha$ is
    a partial homomorphism of the topological space $X$ if the domain $d(\alpha)$
    and the range $r(\alpha)$ are clopen subsets of $X$ and $\alpha$
    is a homeomorphism of $d(\alpha)$ onto $r(\alpha)$.\end{defn}
    \begin{defn} A  partial isometry $v$ in $B_i$ is called a normalising partial
    isometry if $vD_iv^*\subseteq D_i$ and $v^*D_iv\subseteq D_i$
    where $D_i$ denotes the masa in $B_i$. \end{defn} \begin{rem}
    The normaliser $N_{D_i}(B_i)$ is the set of normalising
    partial isometries of $D_i$ in $B_i$. For example if $B_i$ is the upper triangular matrix
    algebra then the normaliser is the set of all upper triangular matrices with
    entries either 0 or of absolute value 1 such that each row or
    column has at most one non-zero entry.\end{rem}
    Although the spectrum
    of limit algebras has been described in detail in literature, we will describe it and then
   illustrate it by working out the spectrum
   for some specific examples.
   Let $B$ denote an AF algebra. Let $D \subseteq B$ be a
   canonical masa and let $X$ denote the Gelfand spectrum. If $x \in
   X$, then there is a decreasing sequence of projections
   $\{p_n\}_{n=1}^\infty$ in $C(X)$ with $\cap_{n=1}^\infty \hat{p_n} =
   \{x\}$; where  $\hat{p_n}$ denotes the spectrum of ${p_n}$ in
   $X$. In other words $\hat{p_n}$ is the image of $p_n$ under the Gelfand map.
   Let $\{e_{ij}^{(n)}\}$ be a set of matrix units of $B$ with respect to $D$.
   Then, $p_n$ can be chosen as a diagonal
   matrix unit in $B_n$. Also note that once you have picked the
   projections $p_n$, you cannot be sure that $p_n \in B_n$, only
   that $p_n \in B_{k_n}$ for some $k_n$. This is just as good.
    If $v$ is a matrix unit in $B$ with $x \in
   \widehat{vv^*}$ then there is an $n \in N$ such that for $n\geq N$
   ,$\{v^*p_nv\}_{n\geq N}$ forms a decreasing set of diagonal
   projections and the intersection $\cap_{n=1}^\infty
   \hat{v^*p_nv}$ is a singleton, say $y$. If $\hat {v}$ denotes the graph of $v$,
   we write $\sigma_v(x)=y
   \Leftrightarrow (x,y)\in \hat{v}$. In this way, $v$ is viewed as
   a partial homeomorphism of $X$, with domain
   $r(v)= \widehat{vv^*}$ and range $d(v)= \widehat{v^*v}$. The orbit of
   $x$ is denoted by $[x]$ with $[x]=\{\sigma_v(x):v$ is a matrix unit of $B$ with $x \in
   \widehat{vv^*}\}$. Each of these orbits is countable. If $A
   \subset B$ is a TAF algebra with $A\cap A^* =D$; we define a
   partial order on each equivalence class in $X$. We call $x\leq y$
   if $\sigma_v(x)=y$ for some matrix unit $v \in A$. This is the
   partial order. This is a total order on each equivalence class
   iff $A$ is strongly maximal (This is discussed in Chapter 4).
    Let $R = \cup \{\hat{v}: $ $v$ is a
   matrix unit of $B$\}. Then $R \subset{X \times X}$. $R$ is topologized
   by letting the compact open sets $\hat {v}$ form a base for
   the topology. If  $P = \cup \{\hat{v}: $ $v$ is a
   matrix unit of $A$\}, then $P \subset R$ is called the
   fundamental relation or spectrum of $A$. More precisely, the
   sets $\hat{v}$ form a base for the topology: They turn out to
   be compact in this topology.
   To summarize; let $E_{ij}^k$
   denote the set of points ($x$,$y$) in $X \times X$
   of the form $(\alpha(y),y)$ where $\alpha$ is the partial
   homeomorphism of $X$ induced by $e_{ij}^k$ and $y$ belongs to the
   domain of $\alpha$.
   Then $P = R(A) = \cup\{E_{ij}^k:e_{ij}^k \in A_k, k=1,2,\cdots\}$ denotes the
   topological binary relation of $A$ with the relative topology.
   The topological binary relation $R(A)$ is the
   spectrum of $A$.
   \bigskip
     \goodbreak
     \subsection{Examples of spectrum for certain limit algebras}
    \begin{enumerate} \item Let us examine the action
    of normalising partial isometries on the maximal ideal space (spectrum) of the masa  in an arbitrary factor
    of the limit algebra with a simple example. Let
    us consider $N_{D_7}(T_7)$. Consider an arbitrary normalising partial
    isometry $v$ of $D_7$ in $T_7$; as mentioned in the remark $v$ has
    entries either 0 or of absolute value 1 such that each row or
    column has at most one non-zero entry. The maximal ideal space $X$ of
    $D_7$ is the set of its minimal diagonal projections. So $X =
    \{e_{11},e_{22},e_{33},e_{44},e_{55},e_{66},e_{77}\} \simeq
    \{1,2,3,4,5,6,7\}$.

    Let $ v =\left[\begin{array}{ccccccc}
      0 & 0 & 1 & 0 & 0 & 0 & 0 \\
      0 & 0 & 0 & 0 & 0 & 0 & 1 \\
      0 & 0 & 0 & 0 & 0 & 0 & 0 \\
      0 & 0 & 0 & 0 & 1 & 0 & 0 \\
      0 & 0 & 0 & 0 & 0 & 0 & 0 \\
      0 & 0 & 0 & 0 & 0 & 0 & 0 \\
      0 & 0 & 0 & 0 & 0 & 0 & 0 \
    \end{array}\right] $

    It is easy to check that $vD_7v^*\subseteq D_7$ and $v^*D_7v\subseteq
    D_7$. Next let us study the action of $vxv^*$ for all $x \in X$ .
     We will observe
    that this action induces a partial map $\alpha_v:S\rightarrow T$ where
    $S\subseteq X$ and $T\subseteq X$ such that $\alpha_v(f)(d) =
    f(vdv^*)$; where $f\in X$, $\alpha_v(f)\in X$ and $d \in D$.\vspace{0.1in}
    Now,
    \begin{itemize}
    \item $ve_{11}v^* = 0$,
    \item $ve_{22}v^* = 0$,
    \item $ve_{33}v^* = e_{11}$,
    \item $ve_{44}v^* = 0$,
    \item $ve_{55}v^* = e_{44}$,
    \item $ve_{66}v^* = 0$,
    \item $ve_{77}v^* = e_{22}$.
    \end{itemize} Let $S = \{e_{33},e_{55},e_{77}\}\simeq \{3,5,7\}$ and  $T =
    \{e_{11},e_{44},e_{22}\}\simeq \{1,4,2\}$. Then for  $f\in X$ such that $\alpha_v(f)\in
    X$ there is a partial map $\alpha_v:S\rightarrow T$ where
    $S\subseteq X$ and $T\subseteq X$ such that $\alpha_v(f)(d) =
    f(vdv^*)$, $d \in D$ is given by
     \begin{itemize}
    \item $\alpha_v(f)(e_{33}) = f(e_{11})$,
    \item $\alpha_v(f)(e_{55}) = f(e_{44})$,
    \item $\alpha_v(f)(e_{77}) = f(e_{22})$.
    \end{itemize} Since $S\subseteq X$ and $T\subseteq X$ are clopen subsets of $X$,
    the partial map $\alpha_v:S\rightarrow T$ is a partial homeomorphism. Thus here the
    normalising partial isometry $v$ induces a partial homeomorphism
    on $X$, the maximal ideal space of $D_7$ and this partial
    homeomorphism essentially moves around the elements of the
    maximal ideal space of $D_7$. We do this for all normalising
    partial isometries and consider all the partial homeomorphisms
    induced by them. Evidently the spectrum of $T_7$ coincides
    with the graphs of all these  normalising
    partial isometries. The topology on $X$ is generated by taking each
    graph as an open  subset of the spectrum. For the above example
    since we are in the discrete case, the topology is trivial.
    \begin{rem}
    Next we identify the spectrum or the topological fundamental relation
    for the refinement limit algebra, the standard limit algebra
    and the alternating limit algebra. We will work the identification of
    the spectrum of the refinement limit algebra and the other two will
    then follow easily. This is mentioned in
    ~\cite{scp92bk}.\end{rem}
     \item Let
    $A={\displaystyle\lim_{k}(T_{n_k},\rho_{k})}$ denote the $r^\infty$
    refinement limit algebra, for some positive integer $r$. Let
    $B={\displaystyle\lim_{k}(M_{r^k},\rho_{k})}$. Then as seen
    before $\rho_{k}(e_{ij}^k)=
    \sum_{t=1}^{r}e_{(i-1){r^k}+t,(j-1){r^k}+t}^{k+1}$ where $\{e_{ij}^{k}: 1 \le i,j \le
    r^{k}\}$ is a m.u.s for $M_{r^k}$. Next we see how to index
    the matrix units using multi-indices. Let $[r]= \{1,2,\cdots,r\}$.
    Then if $k=2$, $\underline {i}= (i_1,i_2)$ and $\underline {j}=
    (j_1,j_2)$ are 2-tuples in
    $[r]^2=[r]\times[r]=\{1,2,\cdots,r\}\times\{1,2,\cdots,r\}$. Thus $\{e_{\underline{i},\underline{j}}:
    \underline{i},\underline{j}\in [r]^2 \}$ is a m.u.s for $M_{r^2}$.
    So in general $\{e_{\underline{i},\underline{j}}:
    \underline{i},\underline{j}\in [r]^k \}$ is a m.u.s for
    $M_{r^k}$. Let $D={\displaystyle\lim_{k}(D_{r^k},\rho_{k})}$ be a
    masa which has the m.u.s  $\{e_{\underline{i},\underline{i}}:
    \underline{i} \in [r]^k, k=1,2,\cdots \}$.
    Let $X$ denote the maximal ideal space of $D$. Given $x \in
    X$, there is a unique sequence
    $(e_{\underline{i_1},\underline{i_1}},e_{\underline{i_2},
    \underline{i_2}},e_{\underline{i_3},\underline{i_3}},\cdots)$
    such that $\hat{e_{\underline{i_n},\underline{i_n}}}(x)=1$ for
    all $n$. Conversely each such sequence
    $e_{\underline{i_1},\underline{i_1}}> e_{\underline{i_2},
    \underline{i_2}}> e_{\underline{i_3},\underline{i_3}}> \cdots$
    corresponds to a unique $x \in X$. Thus it is clear
    that each decreasing sequence of minimal projections $q_x^k
    \in D_{r^k}$ corresponds to a unique point $x \in [r]^\infty$
    under the correspondence
    $q_x^k=e_{(x_1,x_2,\cdots,x_k),(x_1,x_2,\cdots,x_k)}$. In this
    way we identify $X$, with the Cantor
    space $[r]^{\infty} = [r]\times[r]\cdots$ with the product
    topology. Next we identify the spectrum $R(B)$. To
    do this we first specify a relationship between successive
    m.u.s and the natural choice is given by
    $e_{\underline{i},\underline{j}}=\sum_{m=1}^{r_k}e_{(i_1,\cdots,i_k,m),(j_1,\cdots,j_k,m)}$.
    Consequently if  $E_{\underline{i},\underline{i}}$ denotes
    the graph of a partial homeomorphism of the maximal ideal
    space of  $D_{r^k}$ induced by
    $e_{\underline{i},\underline{i}}$ and
    $E_{\underline{i},\underline{j}}$ denotes the same induced by $e_{\underline{i},\underline{j}}$
    then $E_{\underline{i},\underline{i}}=\{(u,u):u \in
    \underline{i}\times[r_{k+1}]\times[r_{k+2}]\cdots\}$ and
    $E_{\underline{i},\underline{j}}=\{(u,v): u=\underline{i}\times{w} ,v=\underline{j}\times{w};\ w\in
   [r_{k+1}]\times[r_{k+2}]\cdots\}$. Thus we identify  the spectrum
   $R(B)$ with the topological equivalence relation on
   the Cantor space $X=[r]^{\infty}$ which consists of all pairs
   $(u,v)$ of points whose tails coincide eventually. As mentioned before we can visualise the topological relation
    $R$ for $e_{\underline{i},\underline{j}}$ as $lying$ $over$
    the square $[0,1]\times[0,1]$ by considering the map $\pi(x)=\sum_{k=1}^{\infty}(x_{k}-1)r^{-k}$.
    For details the reader is referred to ~\cite{scp92bk}. In exactly the same way if
   $n_k=r_1r_2r_3 \cdots r_k$ then we may identify  $R\{\lim_{\rightarrow}(M_{n_k},\rho_{k})\}$
   with an analogous equivalence relation on the Cantor space $[r_1]\times[r_2]\cdots$
    We observe
   that the lexicographic ordering given by
   $\underline{i}\leq\underline{j} \Leftrightarrow
   \underline{i}=\underline{j}$, or $i_m=j_m$, for $1 \leq m < n$, $i_n < j_n$,
   determines $T_{n_k}$ in such a way that the embedding coincides
   with the refinement embedding. So
   $R\{\lim_{\rightarrow}(T_{n_k})\}$=
   $R\{e_{\underline{i},\underline{j}}:\underline{i}\leq\underline{j},k=1,2,\cdots\}$
   =$\{(u,v): u=\underline{i}\times{w} ,v=\underline{j}\times{w},\underline{i} \leq \underline{j} ;\ w\in
   [r_{k+1}]\times[r_{k+2}]\cdots\}$.
   \item It is easy to describe the spectrum of  the standard
   limit algebra on the basis of the previous example. We look at
   $B={\displaystyle\lim_{k}(M_{n_k},\sigma_{k})}$from the previous
   example. Let $A$=
    ${\displaystyle\lim_{k}(T_{n_k},\sigma_{k})}$ denote the $r^\infty$
    standard limit algebra. Then the spectrum
    $R(A)$ can be thought
    of as a subset of $R(B)$ which is determined by the reverse
    lexicographic order in a manner analogous to the previous example.
    Again if $n_k=r_1r_2r_3\cdots r_k$ then we may identify  $R(B)$
    with an analogous equivalence relation on the Cantor space
    $X$ as described in the previous example
    and then obtain $R(A)$ as a subset of $R(B)$ using the reverse
    lexicographic ordering.
    Observe that if $x=(x_1,x_2,\cdots)$ and $y=(y_1,y_2,\cdots)$
    are points in the Cantor space $X=[r_1]\times[r_2]\cdots$ then $x \leq y$ in the
    reverse lexicographic order implies that either $x=y$ or
    $x_k<y_k$ where $x_k$ is the rightmost coordinate of $x$ which
    differs from $y_k$.
    \item The spectrum for the alternating
    limit algebra is obtained by alternating the procedures
    for the standard limit algebra and the refinement limit
    algebra in the previous examples. For details of this refer to
    ~\cite{scp92bk}. \end{enumerate}
    \begin{rem} The topological equivalence
   relation or the spectrum is independent of the m.u.s. Stephen Power ~\cite{scp90b} has
   used this fact to prove that the spectrum of an AF C*-algebra
   is a complete isomorphism invariant up to isometric
   isomorphism. But Donsig, Katsoulis and Hudson ~\cite{dhk01} have shown that
   isometric isomorphism is equal to algebraic isomorphism.
    For TAF algebras, since there
   is only one canonical masa, the spectrum is an invariant of the
   algebra. But for a general limit algebra (a limit of digraph algebras),
   the definition of spectrum depends on the choice of a canonical
   masa. Since it is not known in general if any 2 canonical
   masas are inner conjugate, the spectrum may not be independent
   of the choice of masa.
     But if the
   spectra of two AF C*-algebras are isomorphic as topological relations then the algebras
   themselves are isometrically isomorphic.\end{rem}

     %%% ----------------------------------------------------------------------

\section{Ideals in TAF, TUHF and regular canonical subalgebras of
AF C*-algebras}

%%% ----------------------------------------------------------------------
\subsection{Introduction} The study of various types of ideals of
TAF
   algebras has been a topic of recent research. The structure of
   various types of ideals has been studied. Among them are the
   meet-irreducible ideals in ~\cite{dhhls} and ~\cite{mpl97}, join-irreducible ideals in ~\cite{tdh94},
   prime ideals in ~\cite{tdh97}, lie-ideals in ~\cite{hms98}, n-primitive ideals in ~\cite{mpl93},
   Jacobson radical in ~\cite{apd93}. This paper studies the principal ideals and
   gives a large class of limit algebras in which all ideals are principal.
   \bigskip
     \goodbreak
     \subsection{Principal ideals in Digraph algebras}
    We will
   start the study by defining a principal ideal.\begin{defn} An
   ideal in a TAF algebra is called principal if it is generated
   by a single element of the algebra. \end{defn} We will start by
   looking at the algebra of upper triangular matrices, $T_n$.
   The ideal structure in this finite dimensional algebra is
   itself complicated and since we are looking at infinite dimensional
   analogues of these algebras; we will analyze the structure.
   For example, a generic example of an ideal in
   $T_7$ is given by

   $\left[\begin{array}{ccccccc}
     * & * & * & * & * & * & * \\
     0 & * & * & * & * & * & * \\
     0 & 0 & 0 & 0 & 0 & * & * \\
     0 & 0 & 0 & 0 & 0 & * & * \\
     0 & 0 & 0 & 0 & 0 & * & * \\
     0 & 0 & 0 & 0 & 0 & * & * \\
     0 & 0 & 0 & 0 & 0 & 0 & * \
   \end{array}\right]$.

   The above ideal is a subset of matrices that vanish at all
   entries $(i,j)$, for some fixed set in the $7\times 7$ index set
   of $T_7$. The boundary of the
   zero set is described by a certain non-decreasing function on the
   diagonal set $\{1,2,3,\cdots,7\}$. In general each ideal $I$ of
   $T_n$ is described by an order homomorphism
   $\alpha:\{0,1,\cdots,n\}\rightarrow \{0,1,\cdots,n\}$ such
   that $\alpha(k)\leq k$. So $I=I[\alpha] =\{(x_{ij}):x_{ij}=0
   $ whenever $ i>\alpha (j)\}$. This is the space of matrices which
   vanish below the boundary determined by $\alpha$. Next we claim
   that all ideals of $T_n$ are principal. This is a known result,
   but we will include the proof for completeness.
   Before we write a formal proof
   we will make some observations that will guide us towards the
   result. Let $A=(a_{ij})$
   be an element of $T_{n}$; the algebra of upper triangular matrices.
   Let $\{e_{ij}\}$ denote the matrix unit system in $T_n$. Let
   $a_{ij}$ denote a non-zero entry of $A$. Then
   $e_{ii}Ae_{jj}$ yields the matrix with all entries 0 excepting the
   entry $a_{ij}$. Then $e_{ij}(k_1 e_{j,j}+k_2 e_{j,j+1}
   +k_3 e_{j,j+2}+k_4 e_{j,j+3}\cdots + k_l e_{j,n})$;
   ($k_1,k_2,k_3,\cdots,k_l$ denotes constants)
   gives $a_{ij}$ and all the entries in the row containing $a_{ij}$ that are
   to the right of $a_{ij}$.
   Also $(s_1 e_{i,i}+ s_2 e_{i-1,i}+s_3 e_{i-2,i}+s_4 e_{i-3,i}
   \cdots+s_m e_{1,i})e_{ij}$; ($s_1,s_2,s_3,\cdots,s_m$ denotes constants)
   gives $a_{ij}$ and all the entries in the column containing $a_{ij}$ that
   are vertically above $a_{ij}$. So $e_{ij}$ is
   a generator for all the entries in the row containing $a_{ij}$ that are
   to the right of $a_{ij}$ and all the entries in the column containing $a_{ij}$ that
   are
   vertically above $a_{ij}$. In other words, $e_{ij}$ lies in the
   corner of an $L-block$ matrix and is a generator for that
   matrix. The ideals of $T_n$ can be visualized as a combination of
   $L-block$ matrices and so it is apparent what the generator is;
   it is the sum of elements at the corners of the $L-block$
   matrices. Thus every ideal of the algebra is generated by a
   single element of the algebra and so every ideal is a principal
   ideal. Thus $T_n$ is a principal ideal algebra. For a formal
   proof we adopt an approach which gives us foresight in tackling
   the proof in the infinite dimensional analogue.
   In fact we will prove that digraph algebras are principal ideal
   algebras. Since the algebra of upper triangular matrices is a
   subclass of digraph algebras the proof will also work in that
   setting.
     \begin{thm} Let $A$ denote a digraph algebra. If $I$ is an
     ideal in $A$ then $I$ is a principal ideal. \end{thm}
     \textbf{Proof:} Let $X
      =\{1,2,\cdots,n\}\times\{1,2,\cdots,n\}$. The digraph algebra $A$
      is a subalgebra of $M_n$ such that $A$ contains the diagonal matrices. Now given the digraph algebra $A$,
     we can find a subset $P$
   of $X \times X$ such that $P \subseteq X$ and
   $A = \{a \in M_n: a_{ij}=0\Leftrightarrow (i,j) \notin P\}$. We call
   $P$ the support set of $A$. Consider a subset $F$ of $P$
   such that $I =\{a \in A: a_{ij}=0 \Leftrightarrow (i,j) \notin F\}$ is an
   ideal in $A$. We make the following observations:
   \begin{itemize} \item $(i,i) \in P, \forall i$ \item $(i,j) \in
   P$, $(j,k) \in P$ $\Rightarrow$ $(i,k) \in P$ \item
   $(i,j)\in F$, $(j,k)\in P$ $\Rightarrow$ $(i,k) \in F$
   \item $(i,j)\in P$, $(j,k)\in F$ $\Rightarrow$ $(i,k) \in F$
   \end{itemize} To prove that $I$ is a principal ideal we have
   to obtain a single generator for $I$. We claim that $g =
   \sum_{(i,j)\in F}e_{ij}$ is a generator of $I$ where  $e_{ij}$ denotes a
   matrix unit with $(i,j)th$ entry 1 and remaining entries 0. Let $I_g$
   denote the ideal generated by $g$.
   Now if $e_{ij} \in I$ then $e_{ij} = e_{ii}ge_{jj}$; hence $e_{ij} \in I_g$.
   Next consider the sum $\sum_{(i,j)\in F}\alpha_{ij}e_{ij}$. Evidently $\sum_{(i,j)\in
   F}\alpha_{ij}e_{ij}\in I_g$; whence $I_g \subseteq I$. Also in
   the sum each summand is an element of $I$. This implies that
   $I \subseteq I_g$. Thus we have $I = I_g$ and so $g$ is a
   generator of the ideal $I$. This implies that $I$ is a
   principal ideal and so $A$ is a principal ideal algebra.

     \begin{rem} The ideal of $T_7$ shown below can be thought of as a
   combination of 4 $L-block$ matrices with corners at the
   (1,1)th, (2,2)th, (6,6)th and (7,7)th entries respectively

             $\left[\begin{array}{ccccccc}
     * & * & * & * & * & * & * \\
     0 & * & * & * & * & * & * \\
     0 & 0 & 0 & 0 & 0 & * & * \\
     0 & 0 & 0 & 0 & 0 & * & * \\
     0 & 0 & 0 & 0 & 0 & * & * \\
     0 & 0 & 0 & 0 & 0 & * & * \\
     0 & 0 & 0 & 0 & 0 & 0 & * \
   \end{array}\right]$.
  \vspace{0.2in}

  Thus a generator for the above ideal is;

  \vspace{0.2in}
      $\left[\begin{array}{ccccccc}
     1 & 0 & 0 & 0 & 0 & 0 & 0 \\
     0 & 1 & 0 & 0 & 0 & 0 & 0 \\
     0 & 0 & 0 & 0 & 0 & 0 & 0 \\
     0 & 0 & 0 & 0 & 0 & 0 & 0 \\
     0 & 0 & 0 & 0 & 0 & 0 & 0 \\
     0 & 0 & 0 & 0 & 0 & 1 & 0 \\
     0 & 0 & 0 & 0 & 0 & 0 & 1 \
   \end{array}\right]$.

 \end{rem} \vspace{0.2in}
 \bigskip
     \goodbreak
     \subsection{Principal ideals in limit algebras}
   Since the digraph algebra is a principal
   ideal algebra, it is natural to think about ideals in their
   infinite-dimensional analogues, the regular canonical
   subalgebras.
   Recall that these limit algebras have the digraph algebras as their building blocks.
   It is not at
   all apparent at present if all ideals are principal, although
   it seems feasible from the previous result. Next we will study
   the structure of ideals in TAF algebras. It is apparent that UHF
   algebras have no nontrivial ideals since they are all simple.
   The ideal structure of AF algebras was analyzed in ~\cite{ob72}. The
   following results about ideals of AF algebras are important.
   \begin{defn}A closed subspace $S$ of an approximately finite
   C*-algebra $B$ is said to be inductive relative to the chain
   $\{B_k\}$ of $B$ if $S$ is the closed union of the spaces
   $S \cap B_k$ , for $k=1,2,\cdots$; i.e. $S = \overline{\cup_k{(S \cap B_k)}}$
   \end{defn}It is very important to note that all ideals of
   an approximately finite C*-algebra are inductive. The inductivity of ideals
   plays a very important role in their analysis because this lets us study ideals
   $I$ in an AF-algebra $B$ by looking at the finite dimensional
   pieces $I \cap B_k$ of $I$. An elegant proof for the inductivity of ideals is
   given in ~\cite{scp92bk}. The essence of the proof is that the injective *-homomorphism
   $B_k/I \rightarrow B_k/(I \cap B_k)$ is an isometry. So if any
   sequence $(b_k)$ with each element
    $b_k \in B_k$ converges to an element  $b$ in
   $I$, then the isometry forces $b$ to lie in $\overline{\cup_k{(I \cap
   B_k)}}$. We will state this result as a lemma.
   \begin{lem} Every closed ideal of an approximately finite
   C*-algebra $B$ is inductive relative to the subalgebra chain
   $\{B_k\}$.\end{lem}
         TAF algebras and TUHF algebras have a rich ideal structure.
   We have already analyzed the structure of ideals in the finite
   dimensional factors (upper triangular matrices) of TUHF
   algebras. TUHF nest algebras have a particularly nice
   structure. This is also indicated in the next lemma which holds for
   TUHF nest algebras but does not hold in general. Although the next result
   is not used in the main result it exposes the restrictions in dealing with a proof for
   the main result and the subsequent approach. The result is stated
   as a lemma. \vspace{0.1in} \begin{lem}If $I_n$ is an ideal in some finite
   dimensional factor $A_n$ of a TUHF nest algebra $A$, for each $n$,
   and if $I$ is
   the ideal in $A$ generated by the $I_n$, then $I\cap A_n$ is exactly
   $I_n$.\end{lem}
   This tells that the process of generating ideals in TUHF nest
   algebras does not $add$ additional elements in the factors.
   But the above result is not true in general. For the proof of
   the result and an example illuminating that the process of
   generating ideals in algebras other than TUHF nest algebras
   $adds$ additional elements in the factors, the reader is
   referred to ~\cite{tdh94}. This suggests that generators of an ideal
   $I_n$ in one of the finite dimensional factors $A_n$ of a
   limit algebra $A$ may not yield  generators of the
   ideal $I$ generated by $I_n$ in the limit algebra $A$.
   \bigskip
     \goodbreak
     \section{Spectrum of ideals in limit algebras}

       Next we characterize the spectrum for the ideals of TAF
   algebras and TUHF algebras. Let $B$ be an AF C*-algebra and let
   $I$ be a closed two-sided ideal in $B$. By lemma 3.2, $I$ is
   inductive relative to the subalgebra chain $\{B_k\}$ of $B$.
   Let $\{e_{ij}^k\}$  be a matrix unit system for $B$. Next we define the
   spectrum of $I$ to be the subset $R(I)$ of $R(B)$ given by
   $R(I) = \cup\{E_{ij}^k:e_{ij}^k \in I_k, k=1,2,\cdots\}$ with
   the relative topology. Here recall that $E_{ij}^k$ denotes the set of points $(x,y)$ in $X \times X$
   of the form $(\alpha(y),y)$ where $\alpha$ is the partial
   homeomorphism of $X$ induced by $e_{ij}^k$ and $y$ belongs to the
   domain of $\alpha$. The topological binary relation $R(I)$ is the
   spectrum of $I$. To characterize
   the spectrum of $I$ we will use the next lemma which is a version
   of the local spectral
   theorem for bimodules adapted to ideals. The reader is referred to ~\cite{scp92bk}
   for a proof of the local spectral theorem for bimodules.

    \begin{lem} Let $e_{mn}^l$ be an element of the matrix unit
    system $\{e_{ij}^k\}$ associated with the AF C*-algebra $B$.
    If $I$ is an ideal in $B$ and $E_{mn}^l \subseteq R(I)$, then
    $e_{mn}^l\in I$. \end{lem} The next lemma is immediate from
    the inductivity of ideals and the local spectral theorem.
    Again we adapt a version of the spectral theorem for
    bimodules to ideals. The Bimodule spectral theorem is proved
    in ~\cite{scp92bk}.  \begin{lem} Let $I$ and $J$ be ideals in the AF
    C*-algebra $B$. If $R(I)=R(J)$ then $I=J$. \end{lem}

    %%% ----------------------------------------------------------------------

\section{Groupoid terminology for the spectrum}

%%% ----------------------------------------------------------------------
\subsection{Introduction}
    The spectrum $R(B)$ of an AF
    C*-algebra $B$ is an example of an approximately finite
    r-discrete principal groupoid. In this section we will discuss what this
    means. AF C*-algebras are groupoid C*-algebras and it is the structure
     and substructures of the groupoid which helps in the coordinization of these algebras.
     Any closed subalgebra of an AF C*-algebra also has a functional representation in terms of an open subset
     of the spectrum and we will make this observation in this section.
\bigskip
     \goodbreak
     \subsection{Principal Groupoids}
      To begin with let us define a groupoid.\begin{defn} A
    groupoid is a set $G$ together with a subset $G^2 \subset G
    \times G$, a product map $(a,b)\rightarrow ab$ from $G^2$ to
    G, and an inverse map $a \rightarrow a^{-1}$( so that
    $(a^{-1})^{-1}=a$ ) from $G$ to $G$ such that:\begin{itemize} \item If
    $(a,b),(b,c) \in G^2$, then $(ab,c),(a,bc) \in G^2$ and
    $(ab)c=a(bc)$;\item $(b,b^{-1}) \in G^2$ for all $b \in G$, and
    if $(a,b)\in G^2$, then $a^{-1}(ab)=b$,
    $(ab)b^{-1}=a$.\end{itemize} \end{defn} A trivial example of a
    groupoid is a group. Another example of a groupoid is an
    equivalence relation $R$ on a set $X$. In this case $R^2=
    \{((x,y),(y,z)): (x,y),(y,z) \in R\}$. $R$ is a groupoid with the
    product map from $R^2$ to
    $R$ given by $((x,y),(y,z)) \rightarrow (x,z)$ and the inverse
    defined as $(x,y)^{-1} =(y,x)$. Groupoids based on equivalence
    relation are called $principal$ $groupoids$. \begin{defn} The unit space
    $G^0$ of a groupoid $G$ is defined to be the set $\{xx^{-1}:
    x \in G\}$. \end{defn} \begin{defn} The range map is the map
    $r:G \rightarrow G^0$ given by $r(x)= xx^{-1}$ and the source
    map is the map $d:G \rightarrow G^0$ given by $d(x)= x^{-1}x$.
    \end{defn} For principal groupoids the range and the source maps are
    given by \begin{itemize}\item $r(x,y)= (x,y)(x,y)^{-1} = (x,y)(y,x) = (x,x)$
    \item $d(x,y)= (x,y)^{-1}(x,y)= (y,x)(x,y) = (y,y)$
    \end{itemize}
    We say that $G$ is a topological groupoid if $G$ is
    equipped with a suitable topology for which product and
    inversion are continuous . When $G$ is a topological
    groupoid, we also require that the range and the source maps
    are partial homeomorphisms.
    \begin{lem} Let $G$ be a locally compact groupoid. Then each
    of $r$, $d$ is an open map from $G$ onto $G_0$. \end{lem}
    For a proof of the above lemma refer to ~\cite{altp98}.
    For any locally compact groupoid $G$, let $G^{op}$ denote the
    family of open subsets $A$ of $G$ such that the restrictions
    $r_A$, $d_A$ of $r$, $d$ to $A$ are homeomorphisms onto open
    subsets of $G$.
    \begin{defn} An $r-discrete$ groupoid is a locally compact
    groupoid $G$ such that $G^{op}$ is a basis for the topology of
    $G$. \end{defn}
    \begin{defn} G-sets are subsets of a
    topological groupoid $G$
    such that the restrictions of the range and domain functions are
    one-to-one. \end{defn}
    Note that every $A \in G^{op}$ is a G-set.
     \begin{defn}Let $P$ be an open subset of $G$
     containing $G^0$. $P$ is called a partial order in $G$ if $P
     \circ P \subseteq P$ and $P \cap P^{-1} = G^0$. Moreover, if $P \cup P^{-1} =
     G$ then $P$ is called a total order in $G$. If $P \circ P \subseteq P$
     and $P =P^{-1}$, then we call $P$ an equivalence relation on a subgroupoid of $G$
     \end{defn} \begin{lem} If $P$ is a total order on $G$ then $P$ is closed. \end{lem}
     \begin{rem} A TAF algebra is strongly maximal if
     and only if $P$ is totally ordered. For proofs refer to ~\cite{ms89}.
     Also observe that in this case $P$ is clopen. \end{rem}
     \bigskip
     \goodbreak
     \subsection{Groupoid structure of AF C*-algebras}
         Next we  will observe that an AF C*-algebra $B$ can be expressed as $C^*(G)$ for a
    suitable principal groupoid $G$; i.e $B \approx C^*(G)$.
    We begin by reviewing the construction of a groupoid
    *-algebra, $C^*(G)$, from a locally compact, r-discrete,
    principal groupoid $G$. For the construction in a more general
    setting, see ~\cite{ms89}.
     Let  $C_c(G)$ denote the family of continuous complex-valued
    functions with compact support. We can make $C_c(G)$ into a
    topological *-algebra by defining for $f, g \in C_c(G)$, $(f
    \ast g)$ and $f^\ast$ by,
    \begin{itemize}
    \item $(f \ast g)(a,b)=\sum_{((a,c),(c,b))\in G^2}f(a,c)g(c,b)$.
    \item $f^\ast (a,b) = \overline{f(b,a)}$.
    \end{itemize}
    With these operations $C_c(G)$
    is a $\ast-algebra$. Let $C^\ast(G)$ denote the completion of
    $C_c(G)$ in a natural norm (as defined in Muhly and Solel's
    paper ~\cite{ms89}). Since we are assuming that the groupoid $G$ is
    r-discrete and principal, $C^\ast(G)$ can be viewed as a
    subspace of continuous functions on $G$.(~\cite{r80}, 4.2)
    We consider the space of continuous functions with compact support on $G^0$ i.e. $C_c(G^0)$
    and we identify the closure of this space in $C^\ast(G)$  by
    $C_0(G^0)$.
     The next result is a consequence of
    the Spectral theorem of Bimodules; for a proof the reader is referred to ~\cite{ms89}. For the
    sake of completeness we state the result. \begin{lem} Suppose $B \subseteq C^*(G)$ is a
    closed $C^*(G^0)-bimodule$. Let $Q(B)=\{(x,y) \in G: b(x,y)=0; \forall b \in B\}$ and $I(Q)=
     \{b \in C^*(G): b=0 $\ on Q$ \}$. Then $B= I(Q(B))$.\end{lem} \begin{rem} Again by the
    Spectral theorem of Bimodules, for any subalgebra $A$ of $B$ such that $A$ contains $C_0(G^0)
    $(here $C_0(G^0)$ is the analogue of the diagonal matrices in the context of AF C*-algebras),
    there is a subset $P$ of $G$ such
    that $A$ consists of all the elements in $C^{*}(G)$ that are supported on $P$. Consequently
    we denote $A$ by $A(P)$. Thus TUHF, TAF and regular canonical
    algebras have representations in the form
    $A(P)$.\end{rem}
    \bigskip
     \goodbreak
     \subsection{Examples of Groupoid representations of AF C*-algebras}

    We will observe the above representations for a few examples. As we have mentioned before,
    the spectrum $R(B)$ of an AF C*-algebra $B$ is an example of an r-discrete principal
    groupoid. Also we have noted that the underlying space for the groupoid
    $R(B)$ is the maximal ideal space $X$ for a canonical masa; and $R(B) \subset X \times X$.
    Let us recall that given an AF C*-algebra $B$ and a m.u.s,
    each matrix unit, $v$, from the m.u.s acts on the diagonal $D$
    of $B$ by conjugation $(v^*Dv \subseteq D)$. Consequently each
    matrix unit $v$ induces a partial homeomorphism of $X$
    (maximal ideal space of $D$) onto itself. We denoted the graph
    of this homeomorphism by $\hat{v}$. Then the graphs of all the
    partial homeomorphisms induced by matrix units (matrix units
    suffice) is the spectrum $R(B)$ of $B$ and is a groupoid, say $G$.
    Thus $G \subseteq X \times X$ and is an equivalence relation.
    We put a topology on $G$ and this topology is the smallest
    topology in which every $\hat{v}$ is an open subset. $G$ as
    defined above is an example of an r-discrete, principal, topological groupoid.
    The graph, $\hat{v}$, of the partial homeomorphism associated
    with a matrix unit(or a normalising partial isometry) has the
    following properties: \begin{itemize} \item  $(x,y_1)$ and
    $(x,y_2)$ $\in v$ $\Rightarrow$ $y_1=y_2$.\item  $(x_1,y)$ and
    $(x_2,y)$ $\in v$ $\Rightarrow$ $x_1=x_2$.\end{itemize} As
    mentioned before a subset of $G$ with these properties is
    called a G-set.
     From the above discussion we observe that
    an r-discrete principal groupoid $G$ such that $G \subset X \times
    X$, where $X$ is a topological space, satisfies the
    following conditions: \begin{itemize} \item $G$ is a locally
    compact Hausdorff space. \item The map $(x,y) \rightarrow
    (y,x)$ from $G$ to $G$, and the product map from $G^2$ to
    G given by $((x,y),(y,z)) \rightarrow (x,z)$ are
    continuous.\item The map $x \rightarrow (x,x)$ from $X$ to $G$
    is a homeomorphism. \item The unit space $G^0$ is an open
    subset of $G$. \end{itemize} \begin{rem} In the second
    condition the product map is a partially defined map on $G^2$
    and hence carries the relative product topology of the
    topological space $X \times X$. The fourth condition is also
    called $(r-discreteness)$ because it implies that for each $x
    \in X$ the set $G \cap \{(x,y):y \in X\}$ is a discrete space
    in the relative topology. \end{rem} The spectrum of an AF
    C*-algebra satisfies all the above conditions. Next we look at
    the $2^{\infty}$ refinement and the standard embedding algebras and
    describe their spectra as a groupoid $G$ where $G \subset X \times X$ and
    $X$ is a topological space. \begin{itemize} \item
    Let $ X = \{(a_n):a_n \in \{  0,1\}\}$. Set $ G =
    \{(a,b)\in\ X \times X: a_n \ne b_n $ for a finite no: of
     $n's$ \}. The set X is a locally compact, second countable, Hausdorff
    space and $G \subseteq X \times X$ is a second countable, locally
    compact, r-discrete principal groupoid. We will identify G with
    the spectrum of the refinement limit algebra $A$. The groupoid operations
    on $G$ are as follows: If $x=(a,b)$ and
    $y=(c,d)$ are in G then $xy = (a,d)$ if $b=c$ and it is undefined
    otherwise. As seen above the range and domain functions for G are defined by
    $r((x,y))= (x,x)$ and $d((x,y))=(y,y)$. In the present situation there are G-sets
    of the form \begin{eqnarray*}E_{(a_1, a_2, a_3, \cdots ,a_n),(b_1,
    b_2,b_3,\cdots,b_n)} = \{(a_1, a_2, a_3,\cdots,a_n, x_{n+1}, x_{n+2},
    x_{n+3},\cdots),\\(b_1, b_2, b_3,\cdots,b_n, x_{n+1}, x_{n+2},
     x_{n+3},\cdots): x_k \in \{0,1\} \forall k \geq
     n+1\} \end{eqnarray*} These sets
    are compact and open and forms a base for the topology of G.
     $G$ is r-discrete because the unit space $\{(x,x):x\in X \}= G^0$ is open.
    Also $G^2$ $\subseteq$ G
    and is the set $\{((a,b),(b,c)): a,b,c \in X\}$. Let $C_c(G)$
    denote the space of all continuous complex-valued functions with
    compact support on $G$.  Recall that for $ f, g \in C_c(G)$, $f \ast g$ and $f^\ast$
    on $G^2$ are given by \begin{enumerate}
    \item $(f \ast g)(a,b)=\sum_{((a,c),(c,b))\in G^2} f(a,c)g(c,b)$
    \item $f^\ast (a,b) = \overline{f(b,a)}$.\end{enumerate} With these operations $C_c(G)$
    is a $\ast-algebra$. $C^\ast(G)$
    denotes the completion of $C_c(G)$ in a natural norm as observed before.
    Also $C^\ast(G)$ can be viewed as a
    subspace of continuous functions on $G$. For any G-set E,
    $\chi_E$, which denotes the characteristic function of E, is a partial
    isometry in $C^\ast(G)$. Consequently any element of $C^\ast(G)$
    can be written as a norm limit of linear combinations of $f
    \ast\chi_E$, where E is a G-set and $f$ $\in$ $C_c(G^0)$.
    Also if $E$ and $F$ are  G-sets, then $\chi_E*\chi_F=\chi_{EF}$.
    Now, given the refinement limit algebra $A$ we can
    find an open subset $P$ of $G$ containing $G^0$ such that $A$
    has a representation as functions supported on $P$. Let
    $P$ = $\{(a,b): (a_1,a_2,......a_N)$ precedes $(b_1,b_2,..... b_N)$ in
    lexicographic order and $a_n = b_n$ for ${n > N}$\}. Then $P$ is an open subset of G
    containing $G^0$. $P$ is a total order and uniquely
    determines $A$. We observe that $P$ is the spectrum of $A$.
    $A$ is a subalgebra of $C^\ast(G)$
    which satisfies the condition that its meet with its adjoint is
    $C_c(G^0)$ and $A$ = $\{f \in C^\ast(G): f(h)=0,\forall h \in
    G\backslash P\}$. Thus by the Spectral
    theorem of Bimodules we have that $ A \approx A(P)$. \item
    By imitating the above argument
    with $P$ =$\{(a,b): (a_1,a_2,......a_N)$ precedes $(b_1,b_2,..... b_N)$ in
    reverse lexicographic order and $a_n = b_n$ for ${n > N}$\} we obtain
    a representation of the standard limit algebra
    $A$ as $A(P)$.\end{itemize} \begin{rem}The above procedure can
    be generalized to TAF, TUHF and regular canonical algebras.
    Let $A$ denote a TAF, TUHF or a regular canonical algebra. Since we are dealing
    with ideals we note that if $I$ is an ideal of $A$ then we can find an open
    subset $F$ of $P$ such that $I$ has a representation as
    functions supported on $F$. $F$ is called the ideal set of $I$.\end{rem}

%%% ----------------------------------------------------------------------

\section{Principal ideals in regular canonical subalgebras of AF
C*-algebras}

%%% ----------------------------------------------------------------------
 \subsection{Introduction}
   We have already seen that a TAF algebra $A$ can be represented in
   the form $A(P)$ for an open set $P$ of the groupoid for the enveloping $C^*-algebra$.
    Moreover $A(P)$ is strongly maximal if
    and only if $P$ is totally ordered and in this case each factor $ {A_n}$
    of $A$ can be represented as a direct sum of upper triangular
    matrices, each upper triangular matrix obtained from its
    corresponding full matrix algebra. Also for strongly maximal
    TAF algebras the embeddings $j_n:\ {A_n}{\rightarrow}\ {A_{n+1}}$
    are *-extendible to their appropriate full matrix algebras. We
    will first prove that ideals in subalgebras of second countable, locally compact,
    r-discrete principal groupoids are principal ideals. This will be the main
    result in this paper. Since the
    spectrum of regular canonical subalgebras is a
    locally compact, second countable, r-discrete, principal
    groupoid it will follow that the regular canonical subalgebra
    is a principal ideal algebra and so are strongly maximal TAF
    algebras and strongly maximal in factors TUHF algebras. We
    will obtain these results as corollaries to the main result.
    It is important to mention that the proofs to the corollaries will use
    regularity of embeddings and this is the main characteristic of the above
    limit algebras.
    \bigskip
     \goodbreak
     \subsection{Principal ideals in subalgebras of certain groupoid C*-algebras}
   Before we start with the main result, let us recall that the
   spectral theorem for Bimodules (lemma 6.3.1) implies that there is a
   one-to-one correspondence between ideals $I$
   of a subalgebra $A$ of a groupoid C*-algebra $G$ and open subsets
   $F$ of $P$ such that P$\circ$F$\circ$P
   $\subseteq$F. This correspondence is given by  $I$ = $\{f \in
   C^\ast(G): f(h)=0, \forall h \in P\backslash F\}$. $F$ is called
   the ideal set of $I$. The next lemma is also relevant. For a
   proof the reader is referred to ~\cite{ms89}. \begin{lem}For each
   partial order $P$ in $G$, $A(P)$ is a norm closed subalgebra of
   $C^{*}(G)$ containing $C_{0}(G^0)$. Conversely, each subalgebra $A$
   of $C^{*}(G)$ containing $C_{0}(G^0)$ is of the form $A(P)$ for a
   unique partial order $P$. The correspondence $P \mapsto A(P)$
   is an inclusion preserving bijection between the collection of
   partial orders in $G$ and norm closed subalgebras of $C^{*}(G)$
   containing $C_{0}(G^0)$.\end{lem}
   Next we state and prove the main result in this paper.

    \begin{thm} Let $G$ denote a second countable, locally compact,
    r-discrete principal groupoid that admits a
    cover by compact open G-sets. Let $A$ denote a subalgebra
    of $C^{*}(G)$ such that $C_0(G^0)\subseteq A$. Then $A$ is a Principal Ideal
    Algebra.
   \end{thm}
   \noindent
    \textbf{Proof:}
    Since $G$ is second countable, $G$ has a countable basis of open sets.
    Thus every compact open G-set is the union of a finite
    number of these open sets. Consequently if $G$ has a basis of
    compact open $G$ sets then it has a countable basis.
    Now, given $A$, let $P$ denote the
    open subset of $G$ containing $G^0$ such that $A =A(P)$.
    Since $G$ is covered by countably many compact open G-sets and
    $A$ is a subalgebra of $G$ such that $C_0(G^0)\subseteq A$,
    $P$ is covered by the G-sets it contains.
    As mentioned above there is a one-to-one correspondence between ideals $I$
    of $A$ and open subsets $F$ of $P$ such that $P\circ F\circ P
    \subseteq F$. This correspondence is given by  $I$ = $\{f \in
    C^\ast(G): f(h)=0, \forall h \in P\backslash F\}$ where $F$
    denotes the ideal set of $I$. Next there exist countably many compact open G-sets
    $K_{i}$ such that $\cup_{i=1}^\infty
    \{K_{i}\} = F$.
     We will write $\cup_{i=1}^\infty
    \{K_{i}\} = F$ as a countable disjoint union. Let
     \begin{itemize}
    \item${E_1}= {K_1}$,
    \item${E_2}= {K_2}\backslash{K_1}={K_2}\cap{K_1}^{c}$,
    \item${E_3}= {K_3}\backslash({K_1}\cup{K_2})={K_3}\cap({K_1}\cup{K_2})^{c}$,

     $\vdots$

    \item ${E_i}= {K_i}\backslash({K_1}\cup{K_2}\cdots\cup{K_{i-1}})=
    {K_i}\cap({K_1}\cup{K_2}\cdots\cup{K_{i-1}})^{c}$,

    $\vdots$
    \end{itemize}
    These ${E_i}$ are countable, disjoint, compact and open. Also $\cup_{i=1}^\infty
    \{E_{i}\} = F$. For each
    ${E_i}$, $\chi_{E_i}$ denotes the characteristic function of ${E_i}\subseteq
    G$.  Since ${E_i}$ is compact and open, $\chi_{E_i} \in C_c(G)\subseteq
    C^{*}(G)$.
     Next consider the sequence $(\chi_{E_i})$. Now since $\cup_{i=1}^\infty
    \{E_{i}\} = F$ we
    claim that the sequence $(\chi_{E_i})$
    generates $I$. To prove this we first observe that any
    arbitrary deleted G-set ${K_i}$ (or made smaller by the deletion process)
    can be obtained as follows.
    Now, ${E_1}= {K_1}$ and so $\chi_{K_1} = \chi_{E_1}$.
    Next, ${E_2}= {K_2}\backslash{K_1}$ and so
    ${K_2}={E_2}\cup({K_1}\cap{K_2})$. But ${E_1}= {K_1}$ and so
    ${K_2}={E_2}\cup({E_1}\cap{K_2})$. This implies
    $\chi_{K_2}=\chi_{E_2}+\chi_{E_1}\chi_{K_2}$. Note that
    $\chi_{E_1} \in I$ and so $\chi_{E_1}\chi_{K_2} \in I$.
    Consequently $\chi_{E_2}+\chi_{E_1}\chi_{K_2} \in I$ and so
    $\chi_{K_2}\in I$.
     Next, ${E_3}={K_3}\backslash({K_1}\cup{K_2})$ and so
     ${K_3}={E_3}\cup[{K_3}\cap({K_1}\cup{K_2})]$.
     But ${E_1}\cup{E_2}={K_1}\cup{K_2}$ and so
     $\chi_{K_3}=\chi_{E_3}+\chi_{K_3}(\chi_{E_1}+\chi_{E_2})$.
     Again by definition of an ideal, $\chi_{K_3}\in I$. In
     general for any integer $i$, ${E_i}=
     {K_i}\backslash({K_1}\cup{K_2}\cdots\cup{K_{i-1}})$ and so
     ${K_i}={E_i}\cup[{K_{i}}\cap({K_1}\cup{K_2}\cup{K_3}\cdots\cup{K_{i-1}})]$.
     But
     ${E_1}\cup{E_2}\cdots\cup{E_i}={K_1}\cup{K_2}\cdots\cup{K_i}$ and
     so
     $\chi_{K_i}=\chi_{E_i}+\chi_{K_i}(\chi_{E_1}+\chi_{E_2}+\chi_{E_3}+\cdots+\chi_{E_{i-1}})$.
     Consequently by definition of an ideal, $\chi_{K_i}\in I$.
     Thus $\chi_{K_i}\in I$, for all $i$. Therefore, the sequence $(\chi_{E_i})$
    generates $I$.
    Next we claim
    that $\sum_{i=1}^\infty \frac{\chi_{E_i}}{2^i}$ is a generator of
    the ideal $I$. Since the partial sums $\sum_{i=1}^n
    \frac{\chi_{E_i}}{2^i}$, $n=1,2,\cdots$ are norm convergent, let
    ${g}$ = $\sum_{i=1}^\infty \frac{\chi_{E_i}}{2^i}$.
    The convergence of partial sums  $\sum_{i=1}^n
    \frac{\chi_{E_i}}{2^i}$, $n=1,2,\cdots$ assures that $g \in I$.
    Let ${I_{g}}$ denote the ideal generated
    by ${g}$ and  ${E_{g}}$ denote the ideal set of
    ${I_{g}}$. We will prove that ${I_{g}}$ =$I$. It is evident
    that ${I_{g}} \subseteq I$. In other words, ideal set of ${I_{g}}$
    $\subseteq$ ideal set of $I$ i.e. ${E_{g}}$
    $\subseteq$ $F$.
     To prove $I \subseteq$ ${I_{g}}$, it would suffice to prove
    that $\chi_{E_j} \in  I$ implies that $\chi_{E_j} \in
    {I_{g}}$, for any arbitrary $j$. Now we have that $\cup_{i=1}^\infty \{E_{i}\} =
    F$. Next consider $\chi_{E_j} \in I$, for some $j$.
    Let $d{(\chi_{E_j})}$ and
    $r{(\chi_{E_j})}$ denote the domain and range projections of
    $\chi_{E_{j}}$, we claim that $r{(\chi_{E_j})}{g}d{(\chi_{E_j})}=
    \frac{(\chi_{E_j})}{2^j}$.

    Now, $r{(\chi_{E_j})}{g}d{(\chi_{E_j})}$

      = $ r{(\chi_{E_j})}{(\sum_{i=1}^\infty
    \frac{\chi_{E_i}}{2^i})}d{(\chi_{E_j})}$

    = $ \sum_{i=1}^\infty r{(\chi_{E_j})} \frac{\chi_{E_i}}{2^i}
    d{(\chi_{E_j})}$

    = $\frac{\chi_{E_j}}{2^j} $.

    But $g \in I_g$, so this shows $\frac{\chi_{E_j}}{2^j}$ $\in$ ${I_{g}}$ , and
    hence $\chi_{E_j} \in {I_{g}}$. Thus  ${E_j}$
    $\subseteq$ ${E_{g}}$, for all $j$. In other words, the ideal set of $I$ $\subseteq$ ideal set of
    ${I_{g}}$  and so ${I_{g}}$ =$I$.

    Thus if $I$ is an ideal of $A$ then it is a principal
    ideal and the subalgebra $A$ of $G$ is a
    principal ideal algebra.
    \bigskip
    \goodbreak
    \subsection{Principal ideals in regular canonical subalgebras of
    AF C*-algebras}

     We will need to define $subordinate$ of a matrix unit before we get into the corollary.
     Let $A$ be a regular canonical subalgebra of an
    AF C*-algebra $B$ and let us  denote the presentation of $A$
    by ${A_1}\stackrel{\varphi_{1}}{\rightarrow} {A_2}\stackrel{\varphi_{2}}
     {\rightarrow}{A_3}\stackrel{\varphi_{3}}{\rightarrow}{
     A_4}\stackrel{\varphi_{4}}{\rightarrow} \cdots {A_n}\stackrel{\varphi_{n}}{\rightarrow}{A_{n+1}}
     \cdots$ with star injections $\varphi_{n}:A_{n}
     \rightarrow A_{n+1}$. Also for $m > n$, define $\varphi_{m,n}:\
    {A_n}{\rightarrow}\ {A_m}$ to be the embedding from ${A_n}$
    to ${A_m}$; i.e.  $\varphi_{m,n}= \varphi_{m-1} \circ \cdots \circ \varphi_n$.
    \begin{defn} If $v \in A_n$ is a partial isometry and $m \geq
    n$, a partial isometry $u \in A_m$ is a subordinate of
    $\varphi_{m,n}(v)$ if $r(u)\varphi_{m,n}(v)d(u)=u$. \end{defn}

      \begin{cor} Let $A$ denote a regular canonical subalgebra of an
    AF C*-algebra $B$. Then $A$ is a principal ideal algebra.
      \end{cor}
      \textbf{Proof:} Let $A={A_1}\stackrel{\varphi_{1}}{\rightarrow} {A_2}\stackrel{\varphi_{2}}
     {\rightarrow}{A_3}\stackrel{\varphi_{3}}{\rightarrow}{
     A_4}\stackrel{\varphi_{4}}{\rightarrow} \cdots {A_n}\stackrel{\varphi_{n}}{\rightarrow}{A_{n+1}}
     \cdots$ denote an injective direct system (presentation) of C*-algebras
     $A_n,n=1,2,\cdots$ with star-extendible injections $\varphi_{n}:A_{n}
     \rightarrow A_{n+1}$. Since the matrix unit system need not be unique we choose a
    matrix unit system such that every matrix unit in $B_n$ is a
    sum of matrix units in $B_{n+1}$. Let ${D_n}$ denote the diagonal
    of ${B_n}$, for each $n$ and $lim_{n\rightarrow\infty}{D_n}=
    D$ denote the canonical masa in $B$. The spectrum of $B$ has
    the structure of an r-discrete principal groupoid, say $G$, where
    $G \subseteq X \times X$ and
    $X$ is the spectrum of $D$.  Also $G$ is a second countable,
    locally compact, admits a cover
    by compact open G-sets. Note that in this case the compact
    open G-sets are supports of matrix units.
    Thus $G$ satisfies the hypothesis
    of the main result. Consequently we can obtain
    $C_c(G)$, the space of all continuous complex-valued functions with compact support on
    $G$ as a $\ast-algebra$ and  $C^*(G)$ as the completion of
    $C_c(G)$ with respect to a suitable norm. Now, given  $A$, let $P$ denote the
    open subset of $G$ containing $G^0$ such that $A =A(P)$.
    Thus $G$ is a second countable, locally compact,
    r-discrete principal groupoid that admits a
    cover by compact open G-sets and $A$ is a subalgebra
    of $G$ such that $C_0(G^0)\subseteq A$. Consequently $A$ is a Principal Ideal
    Algebra by Theorem 7.2.2.
    \begin{rem} The above proof is more concrete than the proof of the theorem
    in the sense that we have a presentation of
    the regular canonical subalgebra $A$ of the AF C*-algebra $B$ with supports of matrix units
    providing the compact open G-sets. In this spirit we describe
    the process of obtaining a disjoint collection of compact open
    G-sets ${E_i}$ such that $\cup_{i=1}^\infty
    \{E_{i}\} = F$, where $F$ denotes the ideal set of an ideal $I$
    in $A$ in the proof of the above corollary. Once this
    collection is obtained the remaining part of the proof is the same
    as in Theorem 7.2.2.
    As mentioned in the corollary $P$ denotes the
    open subset of $G$ containing $G^0$ such that $A =A(P)$.
    We have already chosen the matrix unit system in such a
    way that every matrix unit in  $ {A_n}$ is a sum of matrix units
    in $ A_{n+1}$. At this stage we look at the intersection of
    $I$ with individual factors ${A_n}$; i.e. ${A_n} \cap I $, for
    each $n$. Now since $\varphi_{n}({A_n})\subseteq {A_{n+1}}$ we will have
    overlapping matrix units in ${A_n} \cap I $ and ${A_{n+1}}
    \cap I $; for each $n$ and consequently intersecting compact open G-sets in
    the sequence of algebras, ${A_n} \cap I; n=1,2,\cdots$. For each
    compact open G-set $E_i$ which is a support of a matrix unit, let $\chi_{E_i}$
    denote the characteristic function on $E_i$. Then $\chi_{E_i}$ is
    a matrix unit in the sequence of algebras ${A_n} \cap I;
    n=1,2,\cdots$. We will obtain a collection of G-sets $\{E_{i}\}$
    from this sequence in such a manner that these G-sets are
    supports of matrix units in the sequence of algebras and
    $\cup_{j=1}^\infty \{E_{i}\} = F$. Let us first list the matrix
    units of  ${A_n} \cap I $; $n=1,2,\cdots$ in order starting
    with matrix units of ${A_1} \cap I $. Let $\chi_{K_1}$,
    $\chi_{K_2}$, $\chi_{K_3}$, $\chi_{K_4}$,$\cdots$ denote this
    list or sequence. This sequence consists of all the matrix units
    of the sequence of algebras ${A_n} \cap I $; $n=1,2,\cdots$ in order.
    Now to achieve our aim of obtaining the required collection of
    G-sets $E_i$ in the sequence of algebras ${A_n} \cap I$;
    $n=1,2,\cdots$ we delete the matrix units in ${A_n} \cap I$;
    $n=2,3,4,\cdots$ which are subordinate to a previous matrix unit.
    After this deletion process let $\chi_{E_1}$, $\chi_{E_2}$,
    $\chi_{E_3}$, $\chi_{E_4}$,$\cdots$ be the sequence of the
    remaining matrix units. In this sequence no matrix unit is a
    subordinate of a previous matrix unit and consequently these
    matrix units are nonoverlapping and so their supporting G-sets
    $E_i$ are nonintersecting. We also have that $\cup_{i=1}^\infty
    \{E_{i}\} = F$. We claim that the sequence $(\chi_{E_i})$
    generates $I$. To prove this we first observe that any
    arbitrary deleted G-set ${K_i}$ can be obtained as follows. Let
    ${E_i}$ be a G-set from $\cup_{j=1}^\infty \{E_{i}\}$ such that
    $\chi_{K_i}$ is subordinate to $\chi_{E_i}$. That is to say that
    if r(${K_i}$) and d(${K_i}$) denote the range and source maps on
    the groupoid then we have ${K_i}$ = r(${K_i}$)${E_i}$d(${K_i}$).
    Now r(${K_i}$) and d(${K_i}$) are subsets of ${G_0}$ which is the
    unit space of $G$. So from the equation ${K_i}$ =
    r(${K_i}$)${E_i}$d(${K_i}$) we observe that ${K_i}$ is in $F$.
    Observe that the fact that ${E_i}$ is in $F$ and that $F$ is an ideal
    set, automatically places ${K_i}$ in $F$.
    So we conclude that
    the sequence $(\chi_{E_i})$ generates $I$.
    From this point on the remaining portion of
    proof is as in Theorem 7.2.2.\end{rem}

    \medskip

    \begin{rem} We recall that strongly maximal TAF algebras, strongly maximal
    in factors TUHF algebras and in general, TAF algebras are regular canonical subalgebras and
    so the theorem holds in those settings. \end{rem}

    %%% ----------------------------------------------------------------------

%\bibliographystyle{amsplain}
%\bibliography{Tsrilalbib}

\providecommand{\bysame}{\leavevmode\hbox
to3em{\hrulefill}\thinspace}

\end{document}